\documentclass[a4paper,11pt]{article}

\usepackage{amsmath}
\usepackage{amsfonts}
\usepackage{amssymb}
\usepackage{latexsym}
\usepackage{epsfig}
\usepackage{graphicx}
\usepackage{oldgerm}
\usepackage{theorem}

\def\C{\mathbb{C}}
\def\R{\mathbb{R}}
\def\Rp{\mathbb{R}_{\geq 0}}
\def\N{\mathbb{N}}
\def\b1{{\bf 1}}
\def\bE{{\bf E}}

\def\cC{{\cal C}}
\def\cH{{\cal H}}
\def\cD{{\cal D}}

\def\rO{{\rm O}}

\def\Xint#1{\mathchoice
{\XXint\displaystyle\textstyle{#1}}%
{\XXint\textstyle\scriptstyle{#1}}%
{\XXint\scriptstyle\scriptscriptstyle{#1}}%
{\XXint\scriptscriptstyle\scriptscriptstyle{#1}}%
\!\int}
\def\XXint#1#2#3{{\setbox0=\hbox{$#1{#2#3}{\int}$}
\vcenter{\hbox{$#2#3$}}\kern-.5\wd0}}

\def\dashint{\Xint-}

\theorembodyfont{\itshape}

\newtheorem{thm}{Theorem}[section]
\newtheorem{lem}[thm]{Lemma}

\newtheorem{prop}[thm]{Proposition}

\newcommand{\SSC}[1]{\section{#1}\setcounter{equation}{0}}
\newcommand{\qed}{\hbox{\rule[-2pt]{3pt}{6pt}}}

\begin{document}

\title{\bf
Three-Parametric Marcenko--Pastur Density}
\author{
Taiki Endo
\footnote{
Department of Physics,
Faculty of Science and Engineering,
Chuo University,
Kasuga, Bunkyo-ku, Tokyo 112-8551, Japan;
e-mail: taiki@phys.chuo-u.ac.jp}, \,
Makoto Katori
\footnote{
Department of Physics,
Faculty of Science and Engineering,
Chuo University,
Kasuga, Bunkyo-ku, Tokyo 112-8551, Japan;
e-mail: katori@phys.chuo-u.ac.jp}
}
\date{8 January 2020}
\pagestyle{plain}
\maketitle

\begin{abstract}
The complex Wishart ensemble is the statistical ensemble
of $M \times N$ complex random matrices with $M \geq N$
such that the real and imaginary parts of each element
are given by independent standard normal variables.
The Marcenko--Pastur (MP) density $\rho(x; r), x \geq 0$
describes the distribution for squares of the
singular values of the random matrices in this ensemble
in the scaling limit $N \to \infty$, $M \to \infty$ with
a fixed rectangularity $r=N/M \in (0, 1]$.
The dynamical extension of the squared-singular-value distribution
is realized by the noncolliding squared Bessel process,
and its hydrodynamic limit provides the
two-parametric MP density $\rho(x; r, t)$ with time $t \geq 0$,
whose initial distribution is $\delta(x)$.
Recently, Blaizot, Nowak, and Warcho{\l}
studied the time-dependent complex Wishart ensemble
with an external source and introduced
the three-parametric MP density
$\rho(x; r, t, a)$ by analyzing the
hydrodynamic limit of the process starting from $\delta(x-a), a > 0$.
In the present paper, we give useful expressions
for $\rho(x; r, t, a)$ and perform a systematic study
of dynamic critical phenomena
observed at the critical time $t_{\rm c}(a)=a$ when $r=1$.
The universal behavior in the long-term limit
$t \to \infty$ is also reported.
It is expected that the present system having
the three-parametric MP density provides
a mean-field model for QCD
showing spontaneous chiral symmetry breaking.

\vskip 0.2cm


\end{abstract}
\vspace{3mm}

\SSC
{Introduction and Main Results}
\label{sec:introduction}
\subsection{Marcenko--Pastur law}
\label{thm:MP_law}

Assume that $M, N \in \N :=\{1,2, \dots\}$, $M \geq N$.
Consider $M \times N$ complex random matrices
$K=(K_{jk})$ such that the real and the imaginary
parts of elements are i.i.d. and normally distributed
with mean $\mu=0$ and variance $\sigma^2=1/2$.
The normal distribution with mean $\mu$ and variance $\sigma^2$
is denoted by $N(\mu, \sigma^2)$ and
when a random variable $X$ obeys $N(\mu, \sigma^2)$,
we write it as $X \sim N(\mu, \sigma^2)$.
Then the present setting is described as
\[
\Re K_{jk} \sim N(0, 1/2), \quad \Im K_{jk} \sim N(0, 1/2),
\quad j=1, \dots, M, \quad k=1, \dots, N.
\]
We consider a statistical ensemble of $N \times N$
Hermitian random matrices $L$ defined by
\begin{equation}
L= K^{\dagger} K,
\label{eqn:L1}
\end{equation}
where $K^{\dagger}$ denotes the Hermitian conjugate
of $K$.
This ensemble of random matrices is called
the {\it complex Wishart random-matrix ensemble} or
the {\it chiral Gaussian unitary ensemble} (chGUE)
(see, for instance, \cite{For10}).
We denote the eigenvalues of $L$ as $X_j^N, j=1, \dots, N$,
which are nonnegative, since
$L$ is nonnegative definite by definition;
$X_j^N \in \Rp := \{x \in \R: x \geq 0\}$.
The positive square roots of them,
$\sqrt{X_j^N}, j=1, \dots, N$ are called
{\it singular values} of
random rectangular matrices $K$.
In other words, the eigenvalue distribution
of the Hermitian random matrices $L$
can be regarded as the distribution of
squares of singular values of the
rectangular complex random matrices $K$
in the complex Wishart random-matrix ensemble.

Let $\cC_{\rm c}(\R)$ be the set of all continuous
real-valued function with compact support on $\R$.
We consider the empirical measure defined by
\begin{equation}
\Xi^N(dx)= \frac{1}{N} \sum_{j=1}^N
\delta_{X_j^N /M}(dx),
\quad x \in \Rp,
\label{eqn:Xi1}
\end{equation}
where $\delta_y(dx)$ denotes a Dirac measure
concentrated on $y$ such that
$\int_{\R} f(x) \delta_y(dx)=f(y)$ for all $f \in \cC_{\rm c}(\R)$.
Then we take the double limit
$N \to \infty$, $M \to \infty$ for each fixed value of
the rectangularity
\begin{equation}
r=\lim_{\substack{N \to \infty, \cr M \to \infty}} \frac{N}{M} \in (0, 1].
\label{eqn:scaling1}
\end{equation}
We can prove that in this scaling limit (\ref{eqn:scaling1}),
the empirical measure (\ref{eqn:Xi1}) converges weakly to
a deterministic measure $\rho(x) dx, x \in \Rp$ in the sense that
$\int_{\Rp} f(x) \Xi^N(dx) \to \int_{\Rp} f(x) \rho(x) dx$
as $N \to \infty$ for any $f \in \cC_{\rm c}(\R)$.
Moreover, the probability density $\rho$
in the limit measure has a finite support in $\R$
and it is explicitly given as a function of
the parameter $ r \in (0, 1]$ as \cite{MP67}
\begin{equation}
\rho(x; r)
= \frac{\sqrt{(x-x_{\rm L}(r))(x_{\rm R}(r)-x)}}{2 \pi r x}
\b1_{(x_{\rm L}(r), x_{\rm R}(r))}(x)
\label{eqn:MP1}
\end{equation}
with
\begin{equation}
x_{\rm L}(r) :=(1-\sqrt{r})^2, \quad
x_{\rm R}(r) :=(1+\sqrt{r})^2.
\label{eqn:xL_xR_r}
\end{equation}
Here $\b1_{\Lambda}(x)$, $\Lambda \subset \R$
is an indicator function such that
$\b1_{\Lambda}(x)=1$ if $x \in \Lambda$,
and $\b1_{\Lambda}(x)=0$ otherwise.
This convergence theorem is known as
the {\it Marcenko--Pastur law} for
the Wishart random-matrix ensemble
\cite{MP67,For10,AGZ10} and
we call (\ref{eqn:MP1})
the {\it Marcenko--Pastur (MP) density}
in this paper.

\subsection{Dynamical extension of MP density}
\label{thm:dynamical_MP_density}

A dynamical extension of the eigenvalue distribution
of the Wishart random-matrix ensemble is realized by
the solution $\{ {X}^N_j(t) \in \Rp : t \geq 0, j=1, 2, \dots, N\}$
of the following system of stochastic differential equations (SDEs),
\begin{align}
d {X}^N_j(t) &=2 \sqrt{ {X}^N_j(t)} dB_j(t) + 2 (\nu+1) dt
\nonumber\\
& \quad
+ 4  {X}^N_j(t) \sum_{\substack{1 \leq k \leq N, \cr k \not=j}}
\frac{1}{ {X}^N_j(t)- {X}^N_k(t)} dt,
\quad j=1, 2, \dots, N, \quad t \geq 0,
\label{eqn:W_process}
\end{align}
where $\nu=M-N$ and $B_j(t), t \geq 0$ are
independent one-dimensional standard Brownian motions
starting from $x^N_j \in \Rp$, $j=1, \dots, N$.
We assume that
$0 \leq x^N_1 \leq x^N_2 \leq \cdots \leq x^N_N < \infty$.
This one-parameter family ($\nu > 0$) of
$N$-particle stochastic processes was
called (the eigenvalue process of)
the {\it Wishart process} by Bru \cite{Bru91}.
It is also called the {\it Laguerre process}
or the {\it noncolliding squared Bessel process}
\cite{KO01,KT04,KT11}.

We set $\nu=(1-r) M =(1-r) N/r$, $r \in (0, 1]$, and consider the empirical
measure of the solution of SDEs (\ref{eqn:W_process}),
\begin{equation}
\Xi^N_t(dx) :=\frac{1}{N} \sum_{j=1}^N \delta_{ {X}^N_j(t)/M}(dx),
\quad x \in \Rp, \quad t \geq 0.
\label{eqn:empirical2}
\end{equation}
If the initial empirical measure satisfies some
moment conditions and converges weakly to
a measure, $\Xi^N_0(dx)=(1/N) \sum_{j=1}^N \delta_{x^N_j/M}(dx)
\to \xi(dx)$ in the limit
$N \to \infty$, $M \to \infty$ with
$r=N/M$ fixed in $(0, 1]$,
it is proved that $\Xi^N_t(dx)$ converges weakly
to a time-dependent deterministic measure,
which we denote here as
$\rho_{\xi}(x; r, t), t \geq 0$ \cite{CDG01,BNW13}.
We define the {\it Green's function} (the resolvent)
$G_{\xi}(z; r, t)$ by the Stieltjes transform of $\rho_{\xi}$,
\[
G_{\xi}(z; r, t) := \int_{\R} \frac{\rho_{\xi}(x; r, t)}{z-x} dx,
\quad z \in \C \setminus \R.
\]
Then we can prove that this solves the following
nonlinear partial differential equation (PDE) \cite{CDG01,BNW13},
\begin{equation}
\frac{\partial G_{\xi}}{\partial t}
=-\frac{\partial G_{\xi}}{\partial z}
+ r \left\{
\frac{\partial G_{\xi}}{\partial z}-2 z G_{\xi} \frac{\partial G_{\xi}}{\partial z}
- G_{\xi}^2 \right\},
\quad t \in [0, \infty),
\label{eqn:PDE1}
\end{equation}
under the initial condition,
$G_{\xi}(z; r, 0)=\int_{\R} \xi(dx)/(z-x)$, $z \in \C \setminus \R$.
Once the Green's function $G_{\xi}(z; r, 0), z \in \C \setminus \R$ is
determined, we can obtained the density function
following the Sokhotski-Plemelj theorem,
\begin{equation}
\rho_{\xi}(x; r, t)
=- \Im \left[ \lim_{\varepsilon \to 0}
\frac{1}{\pi} G_{\xi}(x+i \varepsilon; r, t) \right],
\quad i :=\sqrt{-1}.
\label{eqn:SPtheorem}
\end{equation}
We regard (\ref{eqn:PDE1}) as
the analogy of the complex Burgers equations
in the inviscid limit
({\it i.e.}, the (complex) one-dimensional Euler equation),
and we call the limit process the
{\it hydrodynamic limit} of the Wishart process \cite{BNW13,BNW14}.

The simplest case of the dynamical extension of
the MP density (\ref{eqn:MP1}) with (\ref{eqn:xL_xR_r})
is obtained by setting the initial distribution,
\begin{equation}
\xi(dx)= \delta_0(dx) := \delta(x) dx,
\label{eqn:xi_0}
\end{equation}
that is, all particles are concentrated on the origin;
$x^N_j=0, j \in \N$.
By the method of complex characteristics,
Blaizot, Nowak, and Warcho{\l} \cite{BNW13}
showed that the Green's function for (\ref{eqn:xi_0}),
$G_{\delta_0}(z)=G_{\delta_0}(z;r,t)$, is given by the solution
of the equation
\begin{equation}
z=\frac{1}{G_{\delta_0}(z)} + \frac{t}{1-rt G_{\delta_0}(z)},
\quad z \in \C \setminus \R, \quad r \in (0, 1], \quad t \geq 0.
\label{eqn:G_eq_0}
\end{equation}
They solved (\ref{eqn:G_eq_0}) and using
the Sokhotski-Plemelj theorem (\ref{eqn:SPtheorem})
derived the time-dependent extension of (\ref{eqn:MP1}),
\begin{equation}
\rho(x; r, t) := \rho_{\delta_0}(x; r, t)
= \frac{\sqrt{(x-x_{\rm L}(r, t))(x_{\rm R}(r, t)-x)}}{2 \pi r t x}
\b1_{(x_{\rm L}(r, t), x_{\rm R}(r, t))}(x)
\label{eqn:MP_t}
\end{equation}
with
\begin{equation}
x_{\rm L}(r, t) :=(1-\sqrt{r})^2t, \quad
x_{\rm R}(r, t) :=(1+\sqrt{r})^2t,
\quad t \in (0, \infty).
\label{eqn:xL_xR_r_t}
\end{equation}
Since $\rho(x; r, 1)$ is equal to
the original MP density
(\ref{eqn:MP1}), the above provides a dynamical derivation
of the Marcenko--Pastur law.
The dependence on $t$ of $\rho(x; r, t)$ given by
(\ref{eqn:MP_t}) with (\ref{eqn:xL_xR_r_t}) is very simple,
but we regard this as the
{\it two-parametric} MP density
in the present paper.

\subsection{Main results}
\label{sec:main_results}

Blaizot, Nowak, and Warcho{\l} \cite{BNW14} have studied
the hydrodynamic limit of the Wishart eigenvalue-process
starting from one-parameter family ($a \geq 0$)
of the initial distribution,
\[
\xi(dx)= \delta_a(dx) := \delta(x-a) dx.
\]
They showed that the Green's function,
$G(z)= G_{\delta_a}(z; r, t)$, $a \geq 0$,
is obtained by the solution of the equation,
\begin{equation}
  z =\frac{1}{G_{\delta_a}(z)}
  +\frac{t}{1-r t G_{\delta_a}(z)}
  +\frac{a}{(1- r t G_{\delta_a}(z))^2},
  \quad z \in \C \setminus \R, \quad r \in (0,1], \quad t \geq 0, \quad a \geq 0.
\label{eqn:z_t_1}
\end{equation}
They claimed in \cite{BNW14} that
a proper solution of this equation yields
$\rho_{\delta_a}(x; r, t)$ via the Sokhotski-Plemelj theorem
and showed an illustration (Fig.1 in \cite{BNW14})
of the time dependence of
this density function for a special case with $r=1$ and $a=1$.
The explicit formula of $\rho_{\delta_a}(x; r, t)$ was,
however, not given there.
See also Section 8 in \cite{LWZ16} and Section 3 in \cite{FG16}
for implicit expressions of $G_{\delta_a}$.

We write the solution discussed in \cite{BNW14}
as $\rho(x; r, t, a)$ and call it
the {\it three-parametric Marcenko-Pastur (MP) density}.
The purpose of the present paper is to report
useful expressions and detailed analysis of this density function
$\rho(x; r, t, a)$ on $\Rp$ with three parameters
$r \in (0, 1], t \geq 0$ and $a \geq 0$.

The main theorem of the present paper is the following.
\begin{thm}
\label{thm:main1}
Let
\begin{align}
S(x; r, t, a) &:=
4 a x^3 - \{8 a^2+4a(3r+2)t-t^2\} x^2
\nonumber\\
& \quad
+2 [2 a^3-2a^2(5r-2) t + a\{r(6r-1)+1\} t^2 -(r+1)t^3 ] x
\nonumber\\
& \quad
+(r-1)^2 t^2 \{a^2-a(4r-2)t+t^2\}.
\label{eqn:S1}
\end{align}
For $r \in (0, 1], t >0, a \geq 0$, consider the case such that
the cubic equation with respect to $x$,
\begin{equation}
S(x; r, t, a)=0,
\label{eqn:cubiceq1}
\end{equation}
has three real solutions, $x_1 \leq x_2 \leq x_3$,
where $x_j=x_j(r, t, a), j=1,2,3$.
Define
\begin{equation}
x_{\rm L}(r, t, a) := x_2(r, t, a),
\quad
x_{\rm R}(r, t, a) :=x_3(r, t, a).
\label{eqn:x_LR_3MP}
\end{equation}
Put
\begin{align}
g &= g(x; r, t, a)
\nonumber\\
& :=
-2 x^3+ 3 \{(2r+1)t+6a \} x^2
-3 \Big[ (r-1)\{(2r+1)t-3a\} t - \sqrt{- 3 S} \Big] x +2(r-1)^3 t^3,
\label{eqn:g1}
\end{align}
with $S=S(x; r, t, a)$ given by (\ref{eqn:S1}),
and define
\begin{align}
\varphi &= \varphi(x; r, t, a)
\nonumber\\
& := - \frac{2}{3} \{x-(r-1)t\}
-\frac{2^{1/3}}{3}
\frac{x^2+\{3a-(2r+1)t\}x+t^2(r-1)^2}{g^{1/3}}
-\frac{g^{1/3}}{3 \times 2^{1/3}}.
\label{eqn:varphi}
\end{align}
Then the three-parametric MP density is given by
\begin{equation}
\rho(x; r, t, a)
= \frac{\sqrt{(x-f_{\rm L}(x; r, t, a))(f_{\rm R}(x; r, t, a)-x)}}{2 \pi r x t}
\b1_{(x_{\rm L}(r, t, a), x_{\rm R}(r, t, a))}(x)
\label{eqn:3_MP1}
\end{equation}
with
\begin{equation}
f_{\rm L}(x; r, t, a) := \left( \frac{\sqrt{d_-}+\sqrt{d_+}}{2} -\sqrt{d_0} \right)^2,
\quad
f_{\rm R}(x; r, t, a) := \left( \frac{\sqrt{d_-}+\sqrt{d_+}}{2} + \sqrt{d_0} \right)^2,
\label{eqn:fLR}
\end{equation}
where
\begin{align}
& d_- = t-a-2\sqrt{a \varphi}, \quad
d_+ = t-a+2\sqrt{a \varphi},
\nonumber\\
& d_0
= \varphi+x+\frac{t-a}{2} + \frac{1}{2} \sqrt{d_- d_+}.
\label{eqn:d1}
\end{align}
\end{thm}
\vskip 0.5cm
\noindent{\bf Remark 1} \,
The formula (\ref{eqn:3_MP1}) for the present
three-parametric MP density seems to be similar to
the original MP density (\ref{eqn:MP1})
and the two-parametric MP density (\ref{eqn:MP_t}).
We should note, however, that
$f_{\rm L}$ and $f_{\rm R}$ appearing in (\ref{eqn:3_MP1})
are not equal to the endpoints $x_{\rm L}$ and $x_{\rm R}$
of the support of density and they depend on $x$
as shown by (\ref{eqn:fLR}) with (\ref{eqn:g1}), (\ref{eqn:varphi})
and (\ref{eqn:d1}).
We can see that
\[
\varphi(x; r, t, 0) := \lim_{a \to 0} \varphi(x; r, t, a)=-x+(r-1)t,
\]
and hence $d_{\pm} \to t, d_0 \to rt$ as $a \to 0$.
Then as $a \to 0$, $f_{\rm L} \to x_{\rm L}(r; t)$,
$f_{\rm R} \to x_{\rm R}(r; t)$ with (\ref{eqn:xL_xR_r_t});
that is, the dependence of $f_{\rm L}$ and $f_{\rm R}$
on $x$ vanishes only in this limit.
Theorem \ref{thm:main1} states that for general $a> 0$,
the endpoints $x_{\rm L}$ and $x_{\rm R}$ of the support
for the three-parametric MP density are given by
the suitably chosen solutions (\ref{eqn:x_LR_3MP})
of the cubic equation (\ref{eqn:cubiceq1}) with (\ref{eqn:S1})
as proved in Section \ref{sec:support} below.
That is, the formula (\ref{eqn:3_MP1}) is universal, but
the choice of solutions (\ref{eqn:x_LR_3MP}) depends
on the parameters $r, t, a$.
The equation (\ref{eqn:z_t_1}) seems to be
a simple perturbation of (\ref{eqn:G_eq_0}),
but the solution turns out to have rich structures,
by which we can describe dynamic critical phenomena
at time $t=t_{\rm c} :=a$ for $a >0$, when $r=1$, as shown below.
\vskip 0.3cm

From the view point of the original random matrix theory,
Theorem \ref{thm:main1} gives the limit theorem
for the eigenvalue distribution of random matrix $L$
given by (\ref{eqn:L1}) in the scaling limit (\ref{eqn:scaling1}),
in which $M \times N$ rectangular complex random
matrices $K=(K_{jk})$ are distributed as
\begin{align}
\Re K_{kk} &\sim N(\sqrt{Ma}, t/2), \quad k=1, \dots, N,
\nonumber\\
\Re K_{jk} &\sim N(0, t/2), \quad j=1, \dots, M, \quad k=1, \dots, N, \quad j \not=k,
\nonumber\\
\Im K_{jk} &\sim N(0, t/2), \quad j=1, \dots, M, \quad k=1, \dots, N,
\label{eqn:K2}
\end{align}
$t >0, r \in (0, 1]$.
By (\ref{eqn:L1}), (\ref{eqn:K2}) gives
\[
\bE[L_{jk}]=M(a+t) \delta_{jk}, \quad
j, k=1, \dots, N.
\]
When $a > 0$, such an ensemble of random matrices
will be called the
{\it Wishart ensemble with an external source}
or the {\it non-centered Wishart ensemble},
since even at $t=0$,
the diagonal elements of $L$ have positive means,
$\bE[L_{jj}]=Ma >0, j=1, \dots, N$,
\cite{BH17,KMW11,DKRZ12,HK13,For13}.

\begin{figure}
\begin{center}
\includegraphics[width=0.6\textwidth]{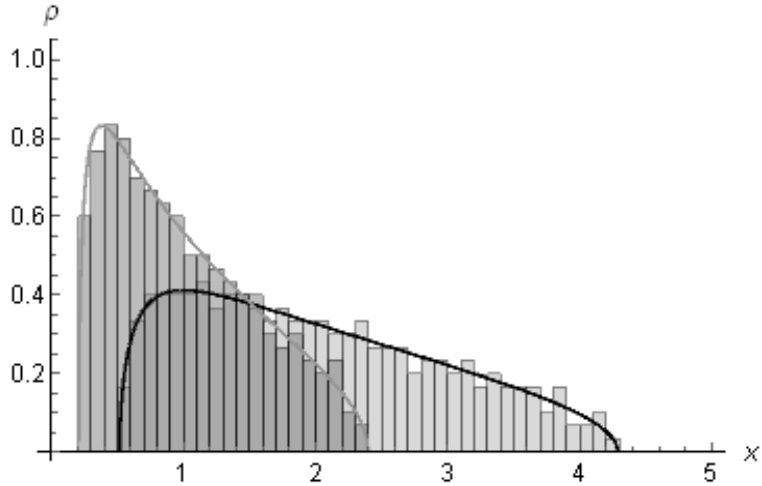}
\end{center}
\caption{
Two histograms for the empirical measures (\ref{eqn:empirical2})
with two different sets of parameters
$(t, a)=(1,0)$ and $(1,1)$
are superposed in order to compare each other.
They show the distributions of
the eigenvalues of $L=K^{\dagger} K$ given by
$K$ of size $1000 \times 300$, whose elements
are randomly generalized following the
probability law (\ref{eqn:K2}) with
$(t, a)=(1, 0)$ and $(1, 1)$.
The original MP density $\rho(x; r=0.3)$ and
the three-parametric MP density
$\rho(x; r=0.3, t=1, a=1)$
are shown by a thin curve and a thick curve, respectively.
Due to an external source at $x=a=1$,
the eigenvalue distribution with $(t, a)=(1,1)$,
which is well fitted by the three-parametric MP density
$\rho(x; r=0.3, t=1, a=1)$, is shifted to the positive direction
and becomes broader compare with the
original MP density $\rho(x; r=0.3)$.
}
\label{fig:histogram}
\end{figure}

In Figure \ref{fig:histogram},
we compare two histograms for
the empirical measures (\ref{eqn:empirical2}) of
the eigenvalues of matrices $L=K^{\dagger} K$
given by $K$ of size $1000 \times 300$
(with the rectangularity $r=300/1000=0.3$),
whose elements are following the
probability law (\ref{eqn:K2}) with different parameters $(t, a)$.
When $(t, a)=(1, 0)$, the distribution of eigenvalues
has a maximum at $x \simeq 0.4$.
When $(t, a)=(1, 1)$, due to an external source at $x=a=1$,
the distribution is shifted to the positive direction
having a maximum at $x \simeq 1$
and becomes broader.
The former is well fitted by
the original MP density $\rho(x; r=0.3)$
and the latter is by the three-parametric MP density
$\rho(x; r=0.3, t=1, a=1)$ given by (\ref{eqn:3_MP1}).

\begin{figure}[htbp]
  \begin{center}
    \begin{tabular}{c}

      \begin{minipage}{0.5\hsize}
        \begin{center}
\includegraphics[width=0.8\textwidth]{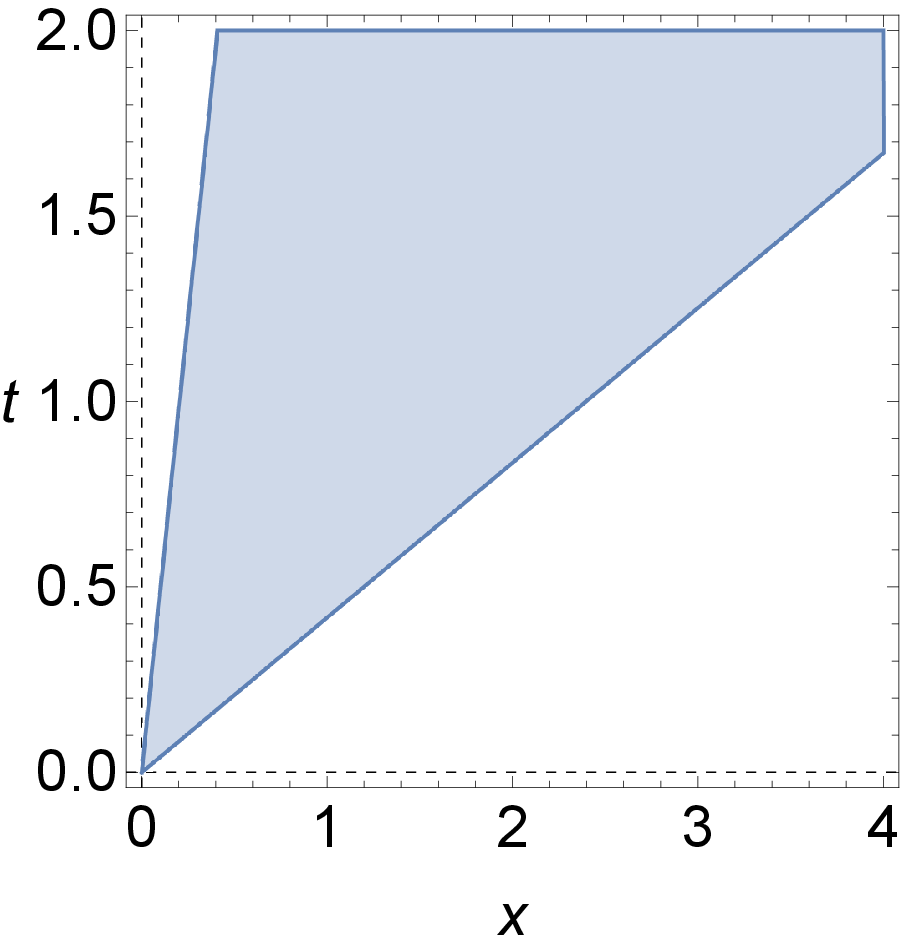}
         \hspace{1.6cm}
        \end{center}
      \end{minipage}

      \begin{minipage}{0.5\hsize}
        \begin{center}
\includegraphics[width=0.8\textwidth]{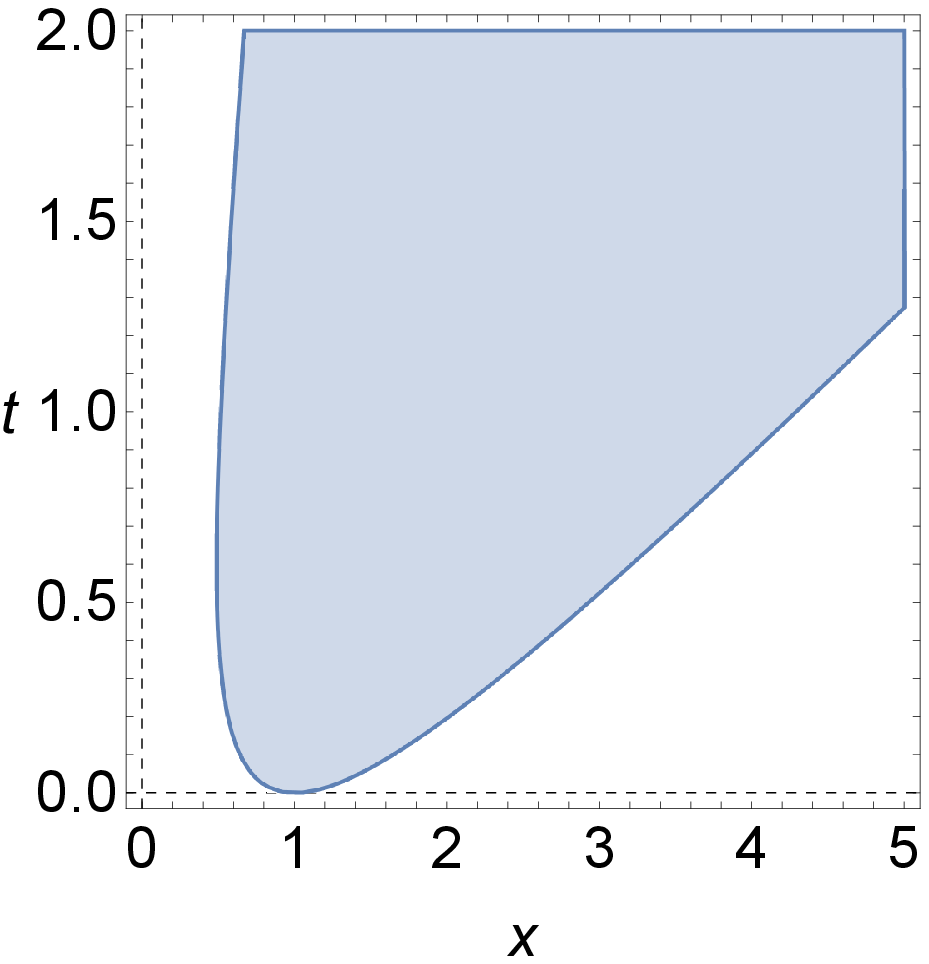}
          \hspace{1.6cm}
        \end{center}
      \end{minipage}

    \end{tabular}
    \caption{
For $r=0.3$, time evolution of the support
$(x_{\rm L}, x_{\rm R})$
is shown on the $(x, t)$-plane
for the two-parametric MP density
$\rho(x; r=0.3, t) := \rho(x; r=0.3, t, a=0)$ in the left,
and for the three-parametric MP density with $a=1$,
$\rho(x; r=0.3, t, a=1)$ in the right.
The supports are extended in time, but
the left edges of supports are kept to be
positive, $x_{\rm L} > 0$, for all $t > 0$.
}
    \label{fig:support_r0.3}
  \end{center}
\end{figure}

\begin{figure}[htbp]
  \begin{center}
    \begin{tabular}{c}

      \begin{minipage}{0.5\hsize}
        \begin{center}
\includegraphics[width=0.8\textwidth]{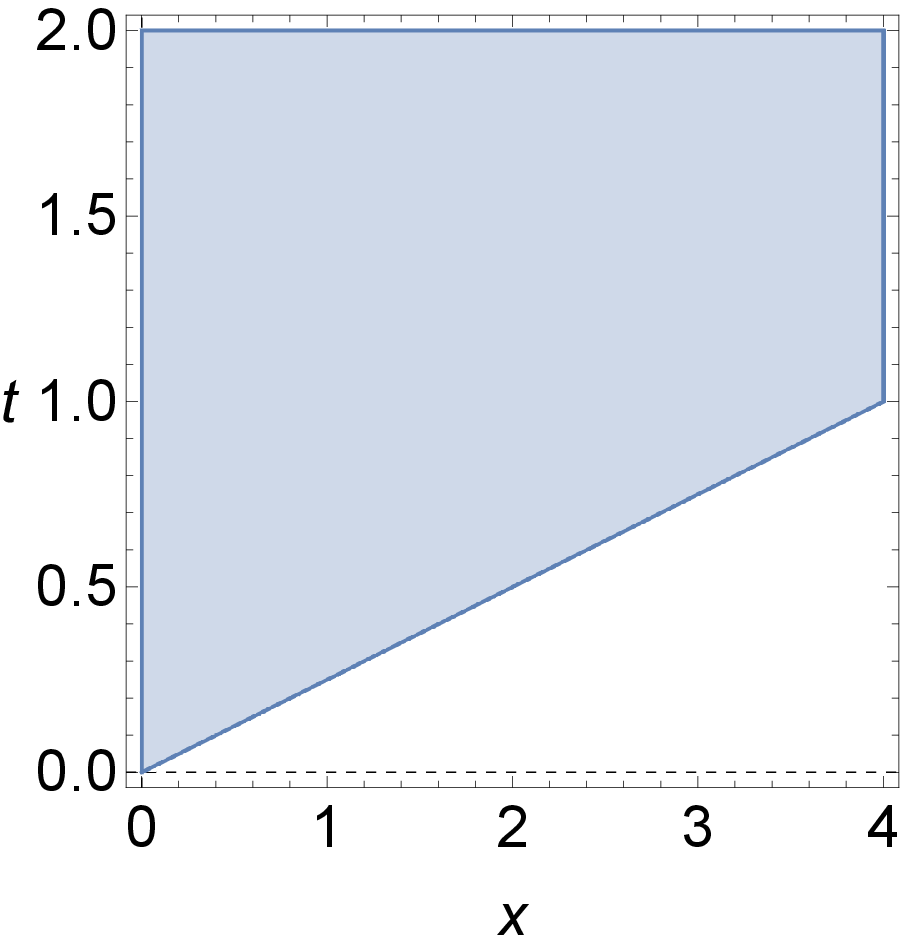}
         \hspace{1.6cm}
        \end{center}
      \end{minipage}

      \begin{minipage}{0.5\hsize}
        \begin{center}
\includegraphics[width=0.8\textwidth]{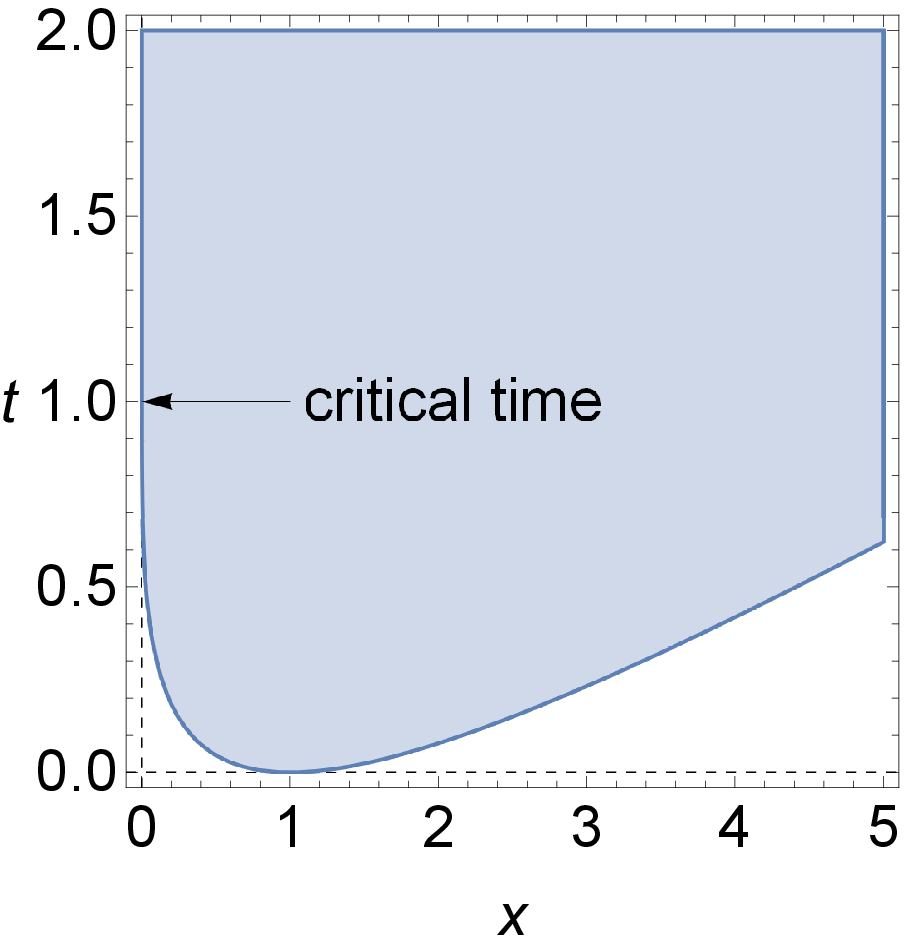}
          \hspace{1.6cm}
        \end{center}
      \end{minipage}

    \end{tabular}
    \caption{
For $r=1$, time evolution of the support
$(x_{\rm L}, x_{\rm R})$
is shown on the $(x, t)$-plane
for the two-parametric MP density
$\rho(x; r=1, t) := \rho(x; r=1, t, a=0)$ in the left,
and for the three-parametric MP density with $a=1$,
$\rho(x; r=1, t, a=1)$ in the right.
In the two-parametric MP density,
the support starts from the singleton $\{ 0 \}$ at $t=0$
and the left edge of support $x_{\rm L}$ is identically zero;
$x_{\rm L} \equiv 1$ for $t \geq 0$.
On the other hand, in the three-parametric MP density with $a=1$,
the support starts from the singleton $\{ 1 \}$ at $t=0$,
and $x_{\rm L} >0$ when $t < t_{\rm c}=1$.
As $t \nearrow t_{\rm c}=1$, however,
$x_{\rm L} \searrow 0$ continuously,
and then $x_{\rm L} \equiv 1$ for $t \geq t_{\rm c}=1$.
We regard $t_{\rm c}=1$ as a critical time.
}
    \label{fig:support_r1}
  \end{center}
\end{figure}

For each values of rectangularity $r \in (0, 1]$
and strength of an external source $a \geq 0$, we can show
time evolution of the support
$(x_{\rm L}(r, t, a), x_{\rm R}(r, t, a) )$ of
$\rho(x; r, t, a)$ on the $(x, t)$-plane, $(\Rp)^2$.
Figure \ref{fig:support_r0.3}
shows the domains
\[
\cD(r, a) := \{ (x_{\rm L}(r, t, a), x_{\rm R}(r, t, a)) : t \geq 0\}
\subset (\Rp)^2
\]
for $(r, a)=(0.3, 0)$ and $(r, a)=(0.3, 1)$.
On the other hand,
Figure \ref{fig:support_r1} shows the domains
for $(r, a)=(1, 0)$ and $(r, a)=(1, 1)$.
As demonstrated by these figures,
we can prove the following qualitative change of the
domain when $a>0$ and $r=1$.


\begin{prop}
\label{thm:support}
Assume that $a > 0$.
\begin{description}
\item{{\rm (i)}} \quad
If and only if $r=1$, $\cD(r,a)$ touches the origin $x=0$.
Otherwise, the left edge of ${\rm supp} \, \rho(x; r, t, a)$
is strictly positive;
$x_{\rm L}(r, t, a) > 0$, $r \in (0, 1)$.
\item{{\rm (ii)}} \quad
When $r=1$, there is a critical time
\[
t_{\rm c}(a)=a
\]
such that
$x_{\rm L}(1, t, a) >0$ while $0 \leq t < t_{\rm c}(a)$, and
$x_{\rm L}(1, t, a) \equiv 0$ for $t \geq t_{\rm c}(a)$.
In particular, just before the
critical time $t_{\rm c}(a)$,
the left edge of ${\rm supp} \, \rho(x; r, t, a)$ behaves as
\[
x_{\rm L}(1, t, a) \simeq \frac{4}{27 a^2}(t_{\rm c}(a)-t)^{\nu}
\quad
\mbox{with $\nu=3$
\quad as $t \nearrow t_{\rm c}(a)$}.
\]
\end{description}
\end{prop}

In the case with $r=1$,
the dynamic critical phenomena
at the critical time $t=t_{\rm c}(a)$
are observed in the vicinity of
the origin as follows.

\begin{figure}
\begin{center}
\includegraphics[width=0.65\textwidth]{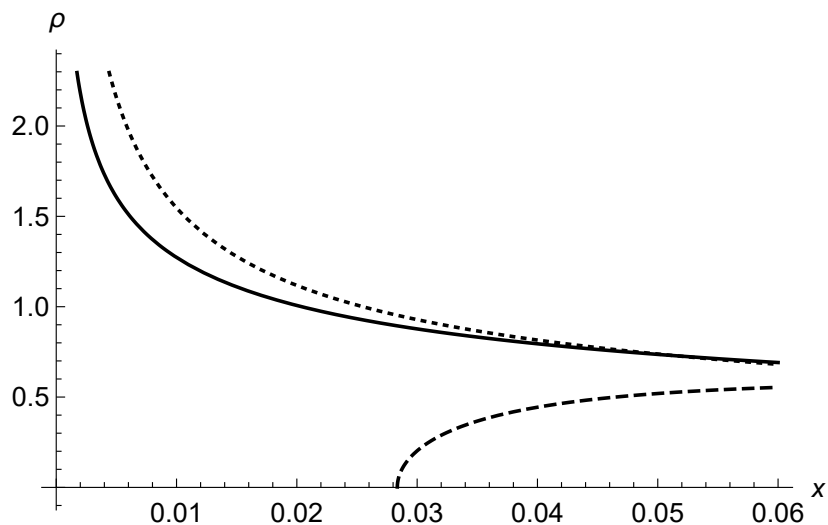}
\end{center}
\caption{
Critical behavior of the three-parametric MP density $\rho$
is shown for $r=1$ and $a=1$ with the critical time $t_{\rm c}(1)=1$.
The dashed curve denotes the emergence of $\rho$
at $x = x_{\rm L} \simeq 0.028$ with the critical exponent
$\beta_1=1/2$ at a subcritical time ($t=0.5 t_{\rm c}(1)$).
The divergence of $\rho$ as $x \searrow 0$ at the
critical time $t=t_{\rm c}(1)$ is shown by a solid curve
and that at a supercritical time ($t=1.5 t_{\rm c}(1)$)
by a dotted curve.
The former with the critical exponent $\gamma_2=1/3$
is weaker than the latter with $\gamma_3=1/2$.
}
\label{fig:rho_critical}
\end{figure}

\begin{prop}
\label{thm:critical}
When $r=1$,
the three-parametric MP density shows
the following dynamic critical phenomena at
$t=t_{\rm c}(a)$.
\begin{description}
\item{{\rm (i)}} \quad
For $0 < t < t_{\rm c}(a)$,
\[
\rho(x; 1, t, a) \simeq C_1(t, a) (x-x_{\rm L}(1, t, a))^{\beta_1}
\quad
\mbox{with $\displaystyle{\beta_1 = \frac{1}{2}}$
\quad as $x \searrow x_{\rm L}(1, t, a)$},
\]
where
\[
C_1(t, a) \simeq \frac{9 a}{4 \pi}(t_{\rm c}(a)-t)^{-\gamma_1}
\quad
\mbox{with $\displaystyle{\gamma_1=\frac{5}{2}}$
\quad as $t \nearrow t_{\rm c}(a)$}.
\]
\item{{\rm (ii)}} \quad
At $t=t_{\rm c}(a)$,
\[
\rho(x; 1, t_{\rm c}(a), a) \simeq \frac{\sqrt{3}}{2 \pi} a^{-2/3} x^{-\gamma_2}
\quad
\mbox{with $\displaystyle{\gamma_2 = \frac{1}{3}}$
\quad as $x \searrow 0$}.
\]
\item{{\rm (iii)}} \quad
For $t > t_{\rm c}(a)$,
\[
\rho(x; 1, t, a) \simeq C_2(t, a) x^{-\gamma_3}
\quad
\mbox{with $\displaystyle{\gamma_3 = \frac{1}{2}}$
\quad as $x \searrow 0$},
\]
where
\[
C_2(t, a) \simeq \frac{1}{\pi t_{\rm c}(a)} (t-t_{\rm c}(a))^{\beta_2}
\quad
\mbox{with $\displaystyle{\beta_2 = \frac{1}{2}}$
\quad as $t \searrow t_{\rm c}(a)$}.
\]
\end{description}
\end{prop}
\vskip 0.5cm
\noindent{\bf Remark 2} \,
The {\it critical exponents} $\nu=3$ and $\gamma_2=1/3$
can be read in the argument given by \cite{BNW14}.
Using the expressions given in Theorem \ref{thm:main1}
here we prove them as well as determining other critical
exponents and critical amplitudes.
The amplitude $C_1$ of $\rho$ in the subcritical time-region
($0<t < t_{\rm c}(a)$) diverges with the critical exponent $\gamma_1=5/2$
as $t \nearrow t_{\rm c}(a)$.
Then both at the critical time ($t=t_{\rm c}(a)$)
and in the supercritical time-region ($t > t_{\rm c}(a)$)
$\rho$ diverges as $x \searrow 0$, but the critical exponents are different.
In the supercritical time-region, the amplitude $C_2$ of the
diverging $\rho$ with the critical exponent $\gamma_3=1/2$
vanishes as $t \searrow t_{\rm c}(a)$ (with the exponent $\beta_2=1/2$).
Consequently the divergence of $\rho$ as $x \searrow 0$
is weakened at $t=t_{\rm c}(a)$
having a smaller value of exponent as
$\gamma_2=1/3 < \gamma_3=1/2$.
See Fig. \ref{fig:rho_critical}.
The proof of Proposition \ref{thm:critical} given in
Subsection \ref{sec:proof3} implies
the scaling relation
\[
\nu=\beta_2+\gamma_1.
\]
\vskip 0.3cm

The functions
$S(x; r, t, a)$, $g(x; r, t, a)$,
$\varphi(x; r, t, a)$, $f_{\rm L}(x; r, t, a)$,
and $f_{\rm R}(x; r, t, a)$,
which appear in Theorem \ref{thm:main1},
are all
{\it homogeneous
as multivariate functions of $x, t, a$}
for each fixed value of $r \in (0, 1]$.
This fact implies the following scaling property
of the three-parametric MP density,
\begin{equation}
\rho(\kappa x; r, \kappa t, \kappa a)
=\frac{1}{\kappa} \rho(x; r, t, a),
\quad r \in (0, 1],
\label{eqn:rho_scaling}
\end{equation}
for an arbitrary parameter $\kappa >0$.
By this property, the following
long-term behavior of the
three-parametric MP density is readily concluded.

\begin{prop}
\label{thm:long_term}
For $a \geq 0$,
\[
\lim_{t \to \infty} \rho(y; r, t, a) dy \Big|_{y=tx}
= \rho(x; r) dx,
\]
where $\rho(x; r)$ is given by
(\ref{eqn:MP1}) with (\ref{eqn:xL_xR_r}).
\end{prop}

\noindent
The long-term behavior
of the present three-parametric MP density
is given by a dilatation of the original MP density
by factor $t$.
In this sense, the original Marcenko--Pastur law
is {\it universal} and it describes the large-scale and long-term
behavior of the Wishart ensemble and process.
\vskip 0.3cm
\begin{figure}
\begin{center}
\includegraphics[width=0.7\textwidth]{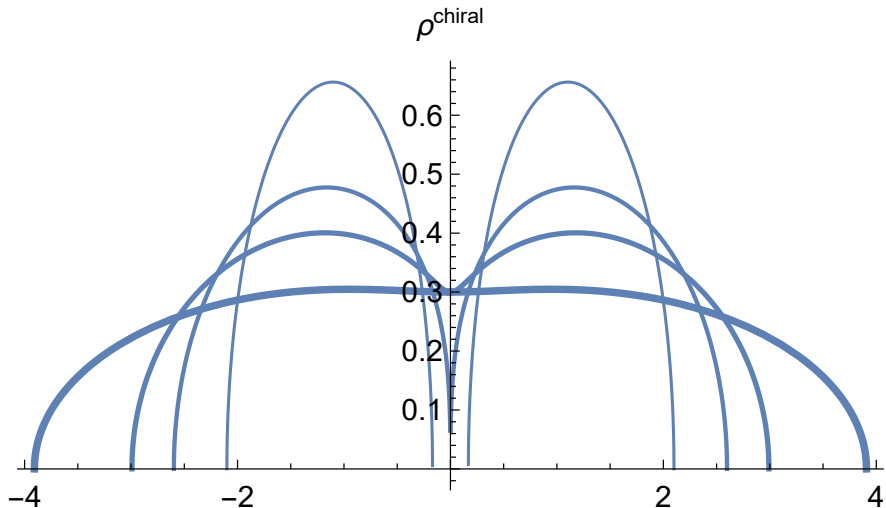}
\end{center}
\caption{
Time evolution of the
hydrodynamical density $\rho^{\rm chiral}(x; r, t, a)$ of the
QCD Dirac operator in the critical case,
which is obtained as (\ref{eqn:rho_chiral})
from the present three-parametric MP density with $r=1$.
For $a=1$, $\rho^{\rm chiral}$ is plotted for
$t=0.5$ (the thinnest curve), 1.0, 1.5
and 3 (the thickest curve).
When $0 \leq t \leq t_{\rm c}=1$, the density at the origin
$\rho^{\rm chiral}(0; r=1, t, a=1) =0$, while
it becomes positive for $t > t_{\rm c}$.
For $t > 3 t_{\rm c}$, $\rho^{\rm chiral}$
shows a relaxation in the sense of
(\ref{eqn:converge}) to the universal density
following Wigner's semicircle law (\ref{eqn:Wigner}).
}
\label{fig:chiral}
\end{figure}
\noindent{\bf Remark 3} \,
So far we have studied the density $\rho$ of
eigenvalues of random matrices $L$ given by (\ref{eqn:L1})
in the hydrodynamic limit.
On the other hand, when the present random matrix ensemble,
chGUE, is applied as a model to
the {\it quantum chromodynamics (QCD)}
in high energy physics,
the density $\rho^{\rm chiral}$ of the
positive-signed and negative-signed singular values of
random rectangular matrix $K$ have been
discussed \cite{JNPZ99,LWZ16,FG16}.
For the transformation from $\rho$ to $\rho^{\rm chiral}$,
see Eq.(3.34) in \cite{FG16}, for instance.
The present three-parametric MP density
$\rho(x; r, t, a)$ given by Theorem \ref{thm:main1}
provides the following hydrodynamical description
of the time-depending spectrum for the
QCD Dirac operator with parameters $r \in (0, 1]$
and $a \geq 0$,
\begin{equation}
\rho^{\rm chiral}(x; r, t, a)
=2 |x| \rho(x^2; r, t, a^2), \quad x \in \R, \quad t > 0,
\label{eqn:rho_chiral}
\end{equation}
under the initial state
\[
\rho^{\rm chiral}(x; r, 0, a)
=\delta(x+a)+\delta(x-a), \quad x \in \R.
\]
Figure \ref{fig:chiral} shows the time evolution of
(\ref{eqn:rho_chiral}) in the critical case $r=1$
with $a=1$.
By (\ref{eqn:rho_chiral}), Proposition \ref{thm:critical} (ii) and (iii) give
the following for $r=1$,
\begin{align}
\rho^{\rm chiral}(x; 1, t_{\rm c}(a), a)
&\simeq \frac{\sqrt{3}}{\pi} a^{-4/3} |x|^{1/\delta}
\quad \mbox{with $\delta=3$ as $|x| \to 0$},
\nonumber\\
\rho^{\rm chiral}(0; 1, t, a)
&\simeq \frac{2}{\pi t_{\rm c}(a^2)}(t-t_{\rm c}(a^2))^{\beta_2}
\quad \mbox{with $\displaystyle{\beta_2=\frac{1}{2}}$
as $t \searrow t_{\rm c}(a)$}.
\label{eqn:MF1}
\end{align}
Moreover, Proposition \ref{thm:long_term} implies
through (\ref{eqn:rho_chiral}) that, when $r=1$,
\begin{align}
\lim_{t \to \infty} \rho^{\rm chiral}(y; 1, t, a) dy
\Big|_{y=\sqrt{t} x}
&= \lim_{t \to \infty} \rho(t x^2; 1, t, a^2) d(t x^2)
\nonumber\\
&= \rho(x^2; 1) d x^2 = 2 x \rho(x^2; 1) dx
\nonumber\\
&= \rho^{\rm Wigner}(x) dx, \quad x \in \R,
\label{eqn:converge}
\end{align}
where
\begin{equation}
\rho^{\rm Wigner}(x)=\frac{1}{\pi} \sqrt{4-x^2}
\b1_{(-2, 2)}(x)
\label{eqn:Wigner}
\end{equation}
is the density function describing
{\it Wigner's semicircle law}
(see, for instance, \cite{For10}).
As mentioned in \cite{LWZ16},
the time evolution of $\rho^{\rm chiral}$ from the
two-peak shape with zero density at the origin
$(0 \leq t \leq t_{\rm c}(a))$
to the universal shape $\rho^{\rm Wigner}$ after
$t \sim 3 t_{\rm c}(a)$ via a critical shape
at $t=t_{\rm c}(a)$ can be interpreted as
a transition from an initial state
with restored chiral symmetry
to a final state with
{\it spontaneous chiral symmetry breaking}.
In \cite{JNPZ99}, we find the argument that
the present system with the density
$\rho^{\rm chiral}(x; r, t, a)$ give
a {\it mean-field model for QCD}
and $\delta=3$ and $\beta_2=1/2$ in (\ref{eqn:MF1}) are the
mean-field values for the scaling exponents
describing a condensation of light quarks
to create massive constituents.
\vskip 0.3cm

The paper is organized as follows.
In Section \ref{sec:proofs} we give proofs
of theorems and propositions given above.
More precisely, Subsections \ref{sec:solving} and \ref{sec:support}
are devoted to the proof of Theorem \ref{thm:main1}.
The proofs of Propositions \ref{thm:support},
\ref{thm:critical}, and \ref{thm:long_term} are given
in Subsections \ref{sec:proof2}, \ref{sec:proof3}, and \ref{sec:proof4},
respectively.
Concluding remarks are given in Section \ref{sec:concluding}.

\SSC
{Proofs of Theorem and Propositions}
\label{sec:proofs}
\subsection{Solving the algebraic equations for
the density and its Hilbert transform}
\label{sec:solving}

For the Green's function
$G_{\delta_a}(z)=G_{\delta_a}(z; r, t)$, $a \geq 0$,
we put
\begin{align*}
  R &=R(x) := \lim_{\varepsilon \to 0} \Re G_{\delta_a}(x+i \varepsilon),
\nonumber\\
  I &=I(x) := - \lim_{\varepsilon \to 0} \Im G_{\delta_a}(x+i \varepsilon),
\end{align*}
that is,
$\lim_{\varepsilon \to 0} G(x+i \varepsilon)=R(x) - i I(x)$.
For the three-parametric MP density
$\rho(x) :=\rho(x; r, t,a)$, its {\it Hilbert transform} is defined by
\begin{align*}
\cH[\rho](x)
&= \frac{1}{\pi} \dashint_{\R}
\frac{\rho(y)}{x-y} dy
\nonumber\\
&:= \frac{1}{\pi} \lim_{\varepsilon \to 0}
\left\{
\int_{-\infty}^{x-\varepsilon} \frac{\rho(y)}{x-y} dy +
\int_{x+\varepsilon}^{\infty} \frac{\rho(y)}{x-y} dy
\right\}.
\end{align*}
The Sokhotski-Plemelj theorem states that
\begin{equation}
\rho(x)=\frac{I(x)}{\pi}, \quad
\cH [\rho(\cdot)](x) = \frac{R(x)}{\pi}.
\label{eqn:SP}
\end{equation}

Let
\begin{equation}
A =\frac{1}{R^2+I^2}, \qquad
B =\frac{1}{(1-rt R)^2+(rt I)^2}.
\label{eqn:AB1}
\end{equation}
By definition, we obtain the equation,
\begin{equation}
A=\frac{(rt)^2}{2 r t R-1 +1/B}.
\label{eqn:AB2}
\end{equation}
For the equation (\ref{eqn:z_t_1}), we obtain the following.
\begin{lem}
\label{thm:equations}
The equation (\ref{eqn:z_t_1}) for the
complex-valued function $G_{\delta_a}(z), a \geq 0$
is equivalent with the following system of
equations for the real-valued functions $A$ and $B$,
\begin{align}
  x & =R A+ (1-rt R)t B + a \left[ (1-rt R)^2-(rt I)^2 \right] B^2,
  \nonumber\\
  0 & =[A-rt^2B]-2a(1-rt R)rt B^2.
\label{eqn:equations1}
\end{align}
\end{lem}
\noindent{\it Proof} \,
We put $z=x+i \varepsilon, x, \varepsilon \in \R$
in (\ref{eqn:z_t_1}) and take the limit $\varepsilon \to 0$.
Then the real part and the imaginary part of
the obtained equation give (\ref{eqn:equations1}).
\qed
\vskip 0.3cm

Before solving the system of equations (\ref{eqn:equations1})
for general $a \geq 0$, first we solve it for
the special case $a=0$.
In this case
(\ref{eqn:equations1}) with (\ref{eqn:SP}) and
(\ref{eqn:AB1}) are simplified as
\begin{align*}
x &= \frac{R_0}{R_0^2+(\pi \rho_0)^2}
+\frac{(1-r t R_0)t}{(1-rtR_0)^2+(rt)^2(\pi \rho_0)^2}
\nonumber\\
&= \left\{ \frac{1}{R_0^2+(\pi \rho_0)^2}- \frac{rt^2}{(1-rtR_0)^2+(rt)^2(\pi \rho_0)^2} \right\} R_0
+\frac{t}{(1-rtR_0)^2+(rt)^2(\pi \rho_0)^2},
\nonumber\\
0 &=\frac{1}{R_0^2+(\pi \rho_0)^2}- \frac{rt^2}{(1-rtR_0)^2+(rt)^2(\pi \rho_0)^2},
\end{align*}
for $R_0 :=R(x; r,t,0)$ and $\rho_0 :=\rho(x; r, t, 0)$.
They give two different expressions for $x$,
\begin{align}
x &= \frac{1}{rt \{R_0^2+(\pi \rho_0)^2\}},
\label{eqn:a0_x1}
\\
x &= \frac{t}{(1-rtR_0)^2+(rt)^2(\pi \rho_0)^2}.
\label{eqn:a0_x2}
\end{align}
From (\ref{eqn:a0_x1}), we have the relation
\begin{equation}
(\pi \rho_0)^2=\frac{1}{r t x}- R_0^2.
\label{eqn:a0_eq1}
\end{equation}
Combining this with (\ref{eqn:a0_x2}), we have
$x=t/(1-2 r t R_0 + rt/x)$,
which is solved as
\begin{equation}
R_0=\frac{x+(r-1)t}{2 rt x}.
\label{eqn:R_0}
\end{equation}
Put (\ref{eqn:R_0}) into (\ref{eqn:a0_eq1}), we obtain
\begin{equation}
(\pi \rho_0)^2
= \frac{1}{(2 rt x)^2}
\{-x^2+2(r+1)tx-(r-1)^2 t^2 \}
= \frac{(x-x_{\rm L})(x_{\rm R}-x)}{(2 r t x)^2},
\label{eqn:sq_2_rho}
\end{equation}
where $x_{\rm L}=x_{\rm L}(r, t)$
and $x_{\rm R}=x_{\rm R}(r, t)$ are
given by (\ref{eqn:xL_xR_r_t})
and the two-parametric MP density
(\ref{eqn:MP_t}) is obtained as
the positive square root of (\ref{eqn:sq_2_rho})
for $x_{\rm L} \leq x \leq x_{\rm R}$.

The above calculation suggests that it will be
easier to obtain $R$ than $I$.
By (\ref{eqn:SP}) and the first equation of (\ref{eqn:AB1}),
if
\[
\frac{1}{A} - R^2 \geq 0,
\]
then
\begin{equation}
\rho(x)=\frac{1}{\pi} \sqrt{ \frac{1}{A}-R^2}.
\label{eqn:rho_Z1}
\end{equation}
Hence if we can express $A$, not using $I$,
but using only $R$ and parameters $r, t, a$,
then the obtained $R$ determines the density
function $\rho$.
Actually we will show that this strategy is successful
in the following.

By eliminating $A$ in (\ref{eqn:equations1}), we
obtain a quadratic equation for $B$ as
\begin{equation}
2a (1-rtR) B^2+(t-a)B-x=0.
\label{eqn:quadratic1}
\end{equation}
We choose the following solution of (\ref{eqn:quadratic1}),
\begin{equation}
B=\frac{a-t+\sqrt{D}}{4a(1-rtR)},
\label{eqn:B1}
\end{equation}
with
\begin{equation}
D=D(x; r, t, a) =8a(1-rtR) x +(t-a)^2,
\label{eqn:Det1}
\end{equation}
by the following reason.
If we put $a=0$, (\ref{eqn:quadratic1}) gives
$B|_{a=0}=x/t$.
On the other hand, (\ref{eqn:Det1}) gives
\[
\sqrt{D} = -(a-t)+4a(1-rtR_0) \frac{x}{t}+\rO(a^2),
\]
and hence (\ref{eqn:B1}) has the correct
limit in $a \to 0$;
$\lim_{a \to 0} B=x/t$.
If we put (\ref{eqn:B1}) into (\ref{eqn:AB2}), then
we have
\begin{equation}
A=\frac{(rt)^2 \{2(2rt R-1)x-a+t-\sqrt{D}\} }{2\{(2rt R-1)^2x+t(2rt R-1)-a\}}.
\label{eqn:A1}
\end{equation}
The function $A$ is indeed expressed
by $R$ and parameters $r, t, a$ apart from $I$.

We find that the second equation of (\ref{eqn:equations1})
gives
\[
tB=\frac{A}{rt}-2a(1-rtR)B^2,
\]
and if we use this equation,
the quadratic equation (\ref{eqn:quadratic1}) for $B$
is written as
\begin{equation}
x=\frac{A}{rt}-aB.
\label{eqn:aAB1}
\end{equation}
Now we put the expression (\ref{eqn:B1})
for $B$ and the expression (\ref{eqn:A1}) for $A$
into (\ref{eqn:aAB1}). Then we obtain the
following equation for $R$,
\begin{align}
  & 8 x (rt R-1)
  \Big\{ 4r^2t^2 xR^2  +2rt^2R -4rt x R + (x -a-t) \Big\}
\nonumber\\
  & \quad \times
  \Big[ 8r^3 t^3 x^2 R^3
  -8r^2t^2 x \{2x+(r-1)t \}R^2 \nonumber \\
  & \quad \qquad +2rt [ 5x^2+ \{(6r-5)t -a\} x+(r-1)^2t^2 ] R \nonumber \\
  & \quad \qquad -[ 2x^2+ \{ (4r -3)t -2a \} x +(r-1) t \{(2r-1)t-a\} ] \Big] =0.
\label{eqn:Eq1}
\end{align}

We want to obtain $R=R(x; r, t, a)$ which solves (\ref{eqn:Eq1})
and satisfy the following
{\it continuity condition with respect to $a$},
\begin{equation}
\lim_{a \to 0} R(x; r, t, a)=R_0
\label{eqn:R_02}
\end{equation}
with (\ref{eqn:R_0}).
This is given as a real solution of the cubic equation
obtained from the last factor in (\ref{eqn:Eq1}),
\begin{align}
  & 8r^3 t^3 x^2 R^3
  -8r^2t^2 x \{2x+(r-1)t \}R^2 \nonumber \\
  & \quad \qquad +2rt [ 5x^2+ \{(6r -5)t -a\} x+(r-1)^2t^2 ] R \nonumber \\
  & \quad \qquad -[ 2x^2+ \{(4r -3)t -2a\} x +(r-1) t \{(2r-1)t-a\} ] =0.
\label{eqn:Eq2}
\end{align}
Applying the Cardano formula, we obtain the solution as
\begin{equation}
R(x; r, t, a)
=\frac{2x+(r-1)t}{3rt x}
  - \frac{x^2 +\{ 3a -(2r+1)t \} x + (r-1)^2 t^2}
  {3 \times 2^{2/3}r t g^{1/3} x}
-\frac{g^{1/3}}{6 \times 2^{1/3}rt x},
\label{eqn:R_solution}
\end{equation}
where $g=g(x; r, t, a)$ is given by
(\ref{eqn:g1}) with (\ref{eqn:S1}).

\begin{lem}
\label{thm:R_to_R0}
The solution (\ref{eqn:R_solution}) of (\ref{eqn:Eq2}) satisfies
the continuity condition (\ref{eqn:R_02}) with (\ref{eqn:R_0}).
\end{lem}
\noindent{\it Proof} \,
If we set $a=0$,
(\ref{eqn:g1}) and (\ref{eqn:S1}) become
\[
g_0 :=g(x; r, t, 0)
= -2 x^3+3(2r+1) t x^2-3 \Big[(r-1)(2r+1)t -\sqrt{-3 S_0} \Big] x + 2(r-1)^3 t^3,
\]
with
\begin{equation}
S_0 := S(x; r, t, 0)=
t^2 \{x^2-2(r+1) t x+(r-1)^2 t^2\}.
\label{eqn:h0}
\end{equation}
In this case, the equality
\begin{equation}
g_0=\frac{1}{4} \left\{ x -(r-1) t + \frac{1}{t} \sqrt{- 3 S_0} \right\}^3
\label{eqn:g0}
\end{equation}
is established.
By putting (\ref{eqn:g0}) with (\ref{eqn:h0}) into
(\ref{eqn:R_solution}) with $a=0$, we can verify
(\ref{eqn:R_02}) with (\ref{eqn:R_0}).
\qed
\vskip 0.3cm

We set
\begin{equation}
\varphi=\varphi(x; r, t, a) :=2 x \{ r t R(x; r, t, a)-1\}.
\label{eqn:varphiA}
\end{equation}
Then it is easy to verify that the expression (\ref{eqn:R_solution})
for $R$ is written as (\ref{eqn:varphi}) for $\varphi$
and that
(\ref{eqn:A1}) gives
\[
\frac{1}{A}-R^2=\frac{1}{(2 r t x)^2}
\Big[ 2(t-a+\sqrt{D})x -\varphi^2 \Big],
\]
with
\begin{equation}
D=-4 a \varphi + (t-a)^2.
\label{eqn:Det2}
\end{equation}
Therefore, if
\begin{equation}
2(t-a+\sqrt{D})x -\varphi^2 \geq 0,
\label{eqn:conditionA}
\end{equation}
then
(\ref{eqn:rho_Z1}) gives
\begin{equation}
\rho(x)=\frac{\sqrt{2(t-a+\sqrt{D})x -\varphi^2}}{2 \pi r t x}.
\label{eqn:rhoA1}
\end{equation}

For the expression (\ref{eqn:3_MP1}) for $\rho(x)$ given in
Theorem \ref{thm:main1}, we perform
the further calculation as follows.
It is easy to verify the equality,
\[
2(t-a+\sqrt{D})=\left(\sqrt{d_-}+\sqrt{d_+}\right)^2,
\]
where $d_{\pm}$ are defined by (\ref{eqn:d1}).
Therefore, we can see that
\begin{align*}
& 2(t-a+\sqrt{D}) x-\varphi^2
\nonumber\\
& \quad
=\left\{ \sqrt{x} \left(\sqrt{d_-}+\sqrt{d_+}\right) -\varphi \right\}
\left\{ \sqrt{x} \left(\sqrt{d_-}+\sqrt{d_+}\right) + \varphi \right\}
\nonumber\\
& \quad
=\left\{\left( \sqrt{x}+\frac{\sqrt{d_-} + \sqrt{d_+}}{2} \right)^2-d_0 \right\}
\left\{d_0 - \left( \sqrt{x}-\frac{\sqrt{d_-} + \sqrt{d_+}}{2} \right)^2 \right\}
\nonumber\\
& \quad
= \left( \sqrt{x}+\frac{\sqrt{d_-} + \sqrt{d_+}}{2} + \sqrt{d_0} \right)
\left( \sqrt{x}+\frac{\sqrt{d_-} + \sqrt{d_+}}{2} - \sqrt{d_0} \right)
\nonumber\\
& \quad
\times
\left( \sqrt{d_0} + \sqrt{x}-\frac{\sqrt{d_-} + \sqrt{d_+}}{2} \right)
\left( \sqrt{d_0} - \sqrt{x}+\frac{\sqrt{d_-} + \sqrt{d_+}}{2} \right)
\nonumber\\
& \quad
=(x-f_{\rm L})(f_{\rm R}-x),
\end{align*}
where $d_0$ is given by (\ref{eqn:d1}), and
$f_{\rm L}=f_{\rm L}(x; r, t, a)$ and
$f_{\rm R}=f_{\rm R}(x; r, t, a)$ are given by (\ref{eqn:fLR}).
Hence (\ref{eqn:rhoA1}) is written as (\ref{eqn:3_MP1}),
provided that the condition (\ref{eqn:conditionA}) is
equivalent with the condition
\begin{equation}
x_{\rm L}(r, t, a) \leq x \leq x_{\rm R}(r, t, a),
\label{eqn:conditionB}
\end{equation}
where $x_{\rm L}(r, t, a)$
and $x_{\rm R}(r, t, a)$
are defined by (\ref{eqn:x_LR_3MP}).

\subsection{Determining the support of density function}
\label{sec:support}

Assume that $t >0, x >0$.
Then if and only if the condition (\ref{eqn:conditionA}) is satisfied,
$\rho$ given by (\ref{eqn:rhoA1}) is positive or zero.
And if and only if $\rho \geq 0$,
its Hilbert transform $R/\pi$ given by the second equation
of (\ref{eqn:SP}) and $\varphi$ defined by (\ref{eqn:varphiA}) are
real valued.
By the explicit expression (\ref{eqn:varphi}) with (\ref{eqn:g1})
for $\varphi$, the following is obvious.
\begin{lem}
\label{thm:rho_real}
If and only if $S(x; r, t, a) \leq 0$,
$\rho(x; r, t, a) \geq 0$.
\end{lem}
\vskip 0.3cm

Now we prove the following.
\begin{lem}
\label{thm:rho_zero}
If $S(x; r, t, a)=0$, then $\rho(x; r, t, a)=0$.
\end{lem}
\noindent{\it Proof} \,
The formula (\ref{eqn:rhoA1}) with (\ref{eqn:Det2}) is
written as
\begin{equation}
\rho(x)=\frac{\sqrt{\varphi F(\varphi)}}{2 \pi r t x \sqrt{2(a-t+\sqrt{D})x + \varphi^2}},
\label{eqn:rho_Z}
\end{equation}
with
\begin{align}
F(\varphi) &:= \frac{1}{\varphi}
\{2(a-t+\sqrt{D})x+\varphi^2 \} \{2(t-a+\sqrt{D})x-\varphi^2\}
\nonumber\\
&= -\varphi^3+4(t-a)x \varphi-16a x^2.
\label{eqn:F1}
\end{align}
Hence we consider the condition of $F(\varphi)=0$.
For $\varphi$ defined by (\ref{eqn:varphiA}),
the cubic equation (\ref{eqn:Eq2})
is written as
\begin{align*}
H(\varphi) &:= \varphi^3 + 2 \{x-(r-1)t \} \varphi^2
+ [ x^2+\{(-2r+3)t-a\} x +(r-1)^2 t^2 ] \varphi
\nonumber\\
& \quad
+t x \{x + (r-1)(a-t) \} =0.
\end{align*}
Therefore,
\[
F(\varphi)=0 \quad \Longleftrightarrow \quad
\widetilde{F}(\varphi) := F(\varphi)+H(\varphi) = 0.
\]
Note that the cubic terms of $\varphi$ are canceled
and $\widetilde{F}(\varphi)$ is reduced to be quadratic in $\varphi$.
We obtain
\begin{align*}
\widetilde{F}(\varphi)
&= 2 \{x-(r-1)t \} \varphi^2
+[ x^2 +\{(-2r+7) t - 5a\} x +(r-1)^2 t^2 ] \varphi
\nonumber\\
& \quad
+ x [(t-16a)x+(r-1)t(a-t)]
\nonumber\\
&=2 \{x-(r-1)t \} (\varphi-\varphi_-)(\varphi-\varphi_+),
\end{align*}
where
\begin{align}
\varphi_- &= \varphi_-(x; r, t, a)=
-\frac{x^2+\{(7-2r)t-5a\} x+(r-1)^2 t^2 + \sqrt{\Delta}}{4\{x-(r-1)t\}},
\nonumber\\
\varphi_+ &= \varphi_+(x; r, t, a)=
-\frac{x^2+\{(7-2r)t-5a\} x+(r-1)^2 t^2 - \sqrt{\Delta}}{4\{x-(r-1)t\}},
\label{eqn:varphi_pm}
\end{align}
with
\begin{align*}
\Delta &=x^4 + 2 \{-(2r-3)t + 59 a \} x^3
+[ \{2r(3r-8)+35\} t^2-2a(58r-33) t + 25 a^2 ] x^2
\nonumber\\
& \quad -2 (r-1)^2 \{(2r-3) t + a \} t^2 x + (r-1)^4 t^4.
\end{align*}
The condition $F(\varphi)=0$ is thus written as
$F(\varphi_+) F(\varphi_-)=0$.
On the other hand, by (\ref{eqn:F1}) and (\ref{eqn:varphi_pm}),
we can show that
\begin{align*}
& F(\varphi_+) F(\varphi_-) = H(\varphi_+) H(\varphi_-)
\nonumber\\
& \quad = -\frac{
\{4 x^3 - \{100 a+(12r+13)t\} x^2 + (r-1) t \{(12r+13) t-25 a\} x
- 4 (r-1)^3 t^3\} x^2 }{8\{x-(r-1)t\}^3}
\nonumber\\
& \qquad \times S(x; r, t, a).
\end{align*}
Note that $x-(r-1)t >0$,
if $x>0$, since $r \in (0, 1], t \geq 0$.
Hence the statement of Lemma is concluded. \qed

\vskip 0.3cm
\noindent
{\it Proof of Theorem \ref{thm:main1}} \,
Assume that $t > 0$.
By definition (\ref{eqn:x_LR_3MP}) of $x_{\rm L}(r, t, a)$ and $x_{\rm R}(r, t, a)$,
and by Lemmas \ref{thm:rho_real} and \ref{thm:rho_zero},
we can conclude that
the condition (\ref{eqn:conditionA}) is equivalent with
the condition (\ref{eqn:conditionB}).
Under the condition
$2(t-a+\sqrt{D})x-\varphi^2>0$
(equivalently,
$x_{\rm L}(r, t, a) < x < x_{\rm R}(r, t, a)$),
if $\varphi F(\varphi) \not=0$,
then $2(a-t+\sqrt{D})x+\varphi^2 \not=0$
by the equality given by the first line of (\ref{eqn:F1}).
Hence $\rho(x)$ given by (\ref{eqn:rho_Z}) is finite for $x>0$.
The proof of Theorem \ref{thm:main1} is thus complete. \qed

\subsection
{Proof of Proposition \ref{thm:support}}
\label{sec:proof2}

The constant term in the cubic function $S(x; r, t, a)$
of $x$ given by (\ref{eqn:S1}) becomes 0 for arbitrary
$t > 0$ and $a \geq 0$, if and only if $r=1$.
This implies Proposition \ref{thm:support} (i).

When $r=1$, the cubic equation (\ref{eqn:cubiceq1})
with (\ref{eqn:S1}) becomes
\[
x \{ 4 a x^2- (8 a^2+20 a t - t^2) x + 4(a-t)^3 \} =0.
\]
For $t >0, a>0$, this equation has three real solutions;
$x=0$ and
\begin{align*}
x =x_{\pm}
&:= \frac{1}{8a} \{
8a^2+20at-t^2
\pm \sqrt{(8a^2+20at-t^2)^2-64a(a-t)^3} \}
\nonumber\\
&= \frac{1}{8a} \{
8a^2+20at-t^2
\pm \sqrt{t} (8a+t)^{3/2} \}.
\end{align*}
When $0 < t < a$,
$8a^2+20at-t^2 > 7a^2+20at>0$,
\[
\sqrt{t} (8a+t)^{3/2}
=\sqrt{(8a^2+20at-t^2)^2-64a(a-t)^3}
< 8a^2+20at-t^2,
\]
and hence $0 < x_- < x_+$.
On the other hand, when $t \geq a$,
\[
\sqrt{t} (8a+t)^{3/2}
=\sqrt{(8a^2+20at-t^2)^2-64a(a-t)^3}
\geq |8a^2+20at-t^2|,
\]
and hence $x_- < 0 < x_+$.
Then we can conclude the following
by the definitions of
$x_{\rm L}(r, t, a)$ and $x_{\rm R}(r, t, a)$
given by (\ref{eqn:x_LR_3MP}).

\begin{lem}
\label{eqn:transition1}
Assume that $a > 0$.
\begin{description}
\item{\rm (i)} \quad
When $0 < t < a$,
\begin{align}
x_{\rm L}(1, t, a) &= \frac{1}{8 a} \{
8 a^2 + 20 a t - t^2 -\sqrt{t} (8a+t)^{3/2} \},
\nonumber\\
x_{\rm R}(1, t, a) &= \frac{1}{8 a} \{
8 a^2 + 20 a t - t^2 +\sqrt{t} (8a+t)^{3/2} \},
\label{eqn:edges_t<a}
\end{align}
and
\[
0 < x_{\rm L}(1, t, a) < x_{\rm R}(1, t, a).
\]

\item{\rm (ii)} \quad
When $t \geq a$,
\begin{align*}
x_{\rm L}(1, t, a) &= 0,
\nonumber\\
x_{\rm R}(1, t, a) &= \frac{1}{8 a} \{
8 a^2 + 20 a t - t^2 +\sqrt{t} (8a+t)^{3/2} \} > 0.
\end{align*}
\end{description}
\end{lem}

Put $t=a-\varepsilon, 0 < \varepsilon \ll 1$
in $x_{\rm L} $ given by (\ref{eqn:edges_t<a}).
Then it is easy to verify that
\begin{equation}
x_{\rm L}(1, t, a)=\frac{4}{27} \frac{\varepsilon^3}{a^2}
+ \frac{8}{81} \frac{\varepsilon^4}{a^3}
+\frac{52}{729} \frac{\varepsilon^5}{a^4}
+\rO(\varepsilon^6).
\label{eqn:xLexp}
\end{equation}
Hence Proposition \ref{thm:support} (ii) is proved.

\subsection
{Proof of Proposition \ref{thm:critical}}
\label{sec:proof3}
In the case with $r=1$, we have the following expressions
from (\ref{eqn:S1}), (\ref{eqn:g1}), (\ref{eqn:varphi}), and (\ref{eqn:Det2}),
\begin{align}
S_1(x) &:= S(x; 1, t, a)
=x [ 4ax^2-(8a^2+20 a t-t^2) x +4(a-t)^3],
\nonumber\\
g_1(x) &:= g(x; 1, t, a)
= x \left[ -2 x^2+9(t+2a)x + 3 \sqrt{-3 S_1(x)} \right],
\nonumber\\
\varphi_1(x) &:= \varphi(x; 1, t, a)
= -\frac{2}{3} x - \frac{2^{1/3}}{3} \frac{\{x+3(a-t)\}x}{g_1(x)^{1/3}}
-\frac{g_1(x)^{1/3}}{3 \times 2^{1/3}},
\nonumber\\
D_1(x) &:= - 4 a \varphi_1(x)+(a-t)^2.
\label{eqn:r1_functions}
\end{align}
Here we write
$\rho_1(x):= \rho(x; 1, t, a)$.

First assume $0 < t < t_{\rm c}(a)=a$.
Put $x_{\rm L} := x_{\rm L}(1, t, a)$ and let $0 < \delta \ll 1$.
Since $S_1(x_L)=0$, we have the expansions in the form,
\begin{align*}
S_1(x_{\rm L}+\delta)
&= c_1 \delta + c_2 \delta^2 + \rO(\delta^3),
\nonumber\\
g_1(x_{\rm L}+\delta)
&= g_1(x_{\rm L})+c_3 \delta^{1/2}+c_4 \delta + \rO(\delta^{3/2}),
\nonumber\\
\varphi_1(x_{\rm L}+\delta)
&= \varphi_1(x_{\rm L})+ c_5 \delta + \rO(\delta^{3/2}),
\end{align*}
where $c_j, j=1, \dots, 5$ are functions of $t, a, x_{\rm L}$,
but independent of $\delta$.
It should be noted that,
in the expansion of $\varphi_1(x_{\rm L}+\delta)$,
the coefficient of term $\delta^{1/2}$ is proportional to
\[
d:= \{x_{\rm L}+3(a-t)\}x_{\rm L}
\left( \frac{2}{g_1(x_{\rm L})} \right)^{1/3}
-\left( \frac{g_1(x_{\rm L})}{2} \right)^{1/3},
\]
where
$g_1(x_{\rm L})=\{-2 x_{\rm L}+9(t+2a)\}x_{\rm L}^{2}$,
and
we can show that $d \propto S_1(x_{\rm L})=0$.
Then, if we note $\rho_1(x_{\rm L})=0$,
(\ref{eqn:rhoA1}) gives
\begin{equation}
\rho_1(x_{\rm L}+\delta)
=\frac{\delta^{1/2}}{\sqrt{2} \pi t x_{\rm L}}
\sqrt{ t-a+\sqrt{D_1(x_{\rm L})}
-\left( \frac{2a x_{\rm L}}{\sqrt{D_1(x_{\rm L})}}
+\varphi_1(x_{\rm L}) \right) c_5} +\rO(\delta).
\label{eqn:rho_expansion}
\end{equation}
By (\ref{eqn:xLexp}) with $\varepsilon := a-t$,
we see that
\begin{align*}
\varphi_1(x_{\rm L}) &= - \frac{4}{9 a} \varepsilon^2 + \rO(\varepsilon^3),
\nonumber\\
D_1(x_{\rm L}) &= \frac{25}{9} \varepsilon^2 + \rO(\varepsilon^3),
\nonumber\\
c_5 &= - \frac{5a}{3} \varepsilon^{-1} + \rO(1).
\end{align*}
Hence by (\ref{eqn:rho_expansion}),
Proposition \ref{thm:critical} (i) is proved.

Next assume $t=t_{\rm c}(a)=a$.
Then (\ref{eqn:r1_functions}) gives
\begin{align*}
S_1(x) &= a(4x-27a) x^2
=-27 a^2 x^2+\rO(x^3),
\nonumber\\
g_1(x) &= 2 \times 3^3 a x^2 + \rO(x^3),
\nonumber\\
\varphi_1(x) &= - a^{1/3} x^{2/3} + \rO(x),
\nonumber\\
D_1(x) &= - 4 a \varphi_1(x) = 4 a^{4/3} x^{2/3} +\rO(x).
\end{align*}
Then (\ref{eqn:rho_expansion}) proves Proposition \ref{thm:critical} (ii).

Finally assume $t > t_{\rm c}(a)=a$.
Then (\ref{eqn:r1_functions}) gives
\begin{align*}
S_1(x) &= - 4 |\varepsilon|^3 x+\rO(x^2),
\nonumber\\
g_1(x) &= 6 \sqrt{3} |\varepsilon|^{3/2} x^{3/2} + \rO(x^2),
\nonumber\\
\varphi_1(x) &= -\frac{t}{|\varepsilon|} x + \rO(x^{3/2}),
\nonumber\\
D_1(x) &= |\varepsilon|^2 + \rO(x),
\end{align*}
where $|\varepsilon|=-\varepsilon=t-a$.
Then (\ref{eqn:rho_expansion}) proves Proposition \ref{thm:critical} (iii).
The proof of Proposition \ref{thm:critical} is hence complete.

\subsection
{Long-Term Scaling and
Proof of Proposition \ref{thm:long_term}}
\label{sec:proof4}

It is obvious that the functions
$S(x; r, t, a)$, $g(x; r, t, a)$,
$\varphi(x; r, t, a)$, $f_{\rm L}(x; r, t, a)$,
and $f_{\rm R}(x; r, t, a)$,
which appeared in Theorem \ref{thm:main1},
are all homogeneous
as multivariate functions of $x, t, a$ for
each fixed value of $r \in (0, 1]$.
Actually, we see that, for an arbitrary parameter $\kappa>0$,
\begin{align*}
& S(\kappa x; r, \kappa t, \kappa a) = \kappa^4 S(x; r, t, a), \quad
\nonumber\\
&
g(\kappa x; r, \kappa t, \kappa a) = \kappa^3 g(x; r, t, a), \quad
\varphi(\kappa x; r, \kappa t, \kappa a) = \kappa \varphi (x; r, t, a),
\nonumber\\
& f_{\rm L}(\kappa x; r, \kappa t, \kappa a) = \kappa f_{\rm L}(x; r, t, a),
\quad
f_{\rm R}(\kappa x; r, \kappa t, \kappa a) = \kappa f_{\rm R}(x; r, t, a).
\end{align*}
Then the scaling property of the
three-parametric MP density (\ref{eqn:rho_scaling})
is concluded.
If we set $\kappa=1/t$, replace $x$ by $t x =: y$,
and take the limit $t \to \infty$ for a fixed $a >0$,
Proposition \ref{thm:long_term} is proved.

\SSC
{Concluding Remarks}
\label{sec:concluding}

In the present paper, we have studied the time-dependent
complex Wishart ensemble of random matrices
with an external source.
Following Blaizot, Nowak, and Warcho{\l} \cite{BNW14},
we have considered the hydrodynamic limit of the
process of the squared-singular-values of random
complex rectangular matrices with
a rectangularity $r \in (0, 1]$.
We solved the algebraic equation (\ref{eqn:z_t_1})
for the Green's function $G_{\xi}$, which is equivalent with
the nonlinear PDE for $G_{\xi}$ (\ref{eqn:PDE1})
under the initial distribution
$\xi(dx)=\delta_a(dx), a > 0$;
a delta measure concentrated at $x=a >0$.
This algebraic equation (\ref{eqn:z_t_1}) was given
and its solution was studied by \cite{BNW14},
but explicit expressions for the density function
has not been available.
In this paper we called the density function of
this system the {\it three-parametric Marcenko--Pastur (MP)
density}, $\rho(x; r, t, a), r \in (0, 1], t >0, a \geq 0$,
and gave useful expressions to $\rho(x; r, t, a)$ (Theorem \ref{thm:main1}).
As an application of the result, the dynamic critical phenomena
were clarified (Propositions \ref{thm:support} and \ref{thm:critical}),
which are observed at the critical time $t_{\rm c}(a)=a$,
if and only if $r=1$ and $a>0$.
There we have introduced six kinds of
{\it critical exponents},
\[
\nu=3, \quad
\beta_1=\beta_2=\frac{1}{2}, \quad
\gamma_1=\frac{5}{2}, \quad
\gamma_2 =\frac{1}{\delta}=\frac{1}{3}, \quad
\gamma_3=\frac{1}{2},
\]
which represent the singularities
of the dynamic critical phenomena.

The present results can be regarded as
{\it macroscopic descriptions}
of the system and the critical phenomena.
{\it Microscopic descriptions} have been also
studied in several papers \cite{KMW11,DKRZ12,HK13,For13}
for the similar systems and
the associated dynamic critical phenomena.
Connection between these two kinds of descriptions \cite{BNW14}
and the universality of such dynamic critical phenomena
will be studied in more detail in the future.
As mentioned in Remark 3 given at the end of
Section \ref{sec:main_results},
in the context of high energy physics,
the present macroscopic description can be regarded
as a mean-field approximation for
more precise theory of QCD
which exhibits spontaneous
chiral symmetry breaking.

As emphasized in \cite{BNW13,BNW14},
the MP density of the Wishart random-matrix ensemble
has been used in a broad range of mathematical sciences,
physics, and information theory
(see the references in \cite{BNW13,BNW14}).
It is expected that the non-centered Wishart ensembles/processes
and the present three-parametric MP density will be
also useful in many applications, where the mean zero condition
cannot be assumed.

\vskip 1cm
\noindent{\bf Acknowledgements} \,
The present authors thank Hiroya Baba for useful discussion
when the present study was started.
They also thank anonymous referees very much,
who suggested them to discuss the present results
in the context of the application of chGUE to QCD.
This work was supported by
the Grant-in-Aid for Scientific Research (C) (No.19K03674),
(B) (No.18H01124), and
(S) (No.16H06338)
of Japan Society for the Promotion of Science.



\end{document}